\documentclass{amsart}
\usepackage{amssymb, amsmath, latexsym, mathabx, dsfont}

\renewcommand{\baselinestretch}{\baselinestretch}
\renewcommand{\baselinestretch}{1.1}
\numberwithin{equation}{section}

\newtheorem{thm}{Theorem}[section]
\newtheorem{lem}[thm]{Lemma}

\theoremstyle{definition}

\theoremstyle{remark}

\newcommand{\z}{{\mathbb Z}}

\begin{document}
\title{A generalization of Gauss' triangular theorem}
\author{Jangwon Ju and Byeong-Kweon Oh }

\address{Department of Mathematical Sciences, Seoul National University, Seoul 08826, Korea}
\email{jjw@snu.ac.kr}
\thanks{This work of the first author was supported by BK21 PLUS SNU Mathematical Sciences Division.}

\address{Department of Mathematical Sciences and Research Institute of Mathematics, Seoul National University, Seoul 08826, Korea}
\email{bkoh@snu.ac.kr}
\thanks{This work of the second author was supported by the National Research Foundation of Korea (NRF-2014R1A1A2056296).}



\begin{abstract}  A quadratic polynomial $\Phi_{a,b,c}(x,y,z)=x(ax+1)+y(by+1)+z(cz+1)$ is called universal if the diophantine equation $\Phi_{a,b,c}(x,y,z)=n$ has an integer solution
$x,y,z$ for any non negative integer $n$. In this article, we show that if $(a,b,c)=(2,2,6), (2,3,5)$ or $(2,3,7)$, then $\Phi_{a,b,c}( x,y,z)$ is universal. These were conjectured  by Sun in 
\cite {Sun}.
\end{abstract}

\maketitle

\section{introduction}
A triangular number is a number represented as dots or pebbles arranged in the shape of an equilateral triangle. More precisely, the $n$-th triangular number is defined by 
$T_n=\frac{n(n+1)}{2}$ for any non negative integer $n$. 

In 1796, Gauss proved that every positive integer can be expressed as a  sum of three triangular numbers, which was first asserted by Fermat in 1638.
This follows from the Gauss-Legendre theorem, which states that every positive integer which is not of the form $4^k(8l+7)$ with  non negative integers $k$ and $l$, is a sum of three squares of integers.  

In general, a ternary sum $aT_x+bT_y+cT_z\ (a,b,c >0)$ of triangular numbers is called {\it universal} if  for any non negative integer $n$, the diophantine equation
$$
aT_x+bT_y+cT_z=n
$$
has  an integer solution $x,y,z$.  In 1862, Liouville generalized the Gauss' triangular theorem by proving that a ternary sum $aT_x+bT_y+cT_z$ of triangular numbers is universal if and only if $(a,b,c)$ is one of  the following triples:
$$
\begin{array}{llllllll}
&(1,1,1), &(1,1,2), &(1,1,4), &(1,1,5), &(1,2,2), &(1,2,3), &(1,2,4).
\end{array}
$$

Recently,  Sun in \cite{Sun} gave an another generalization of the Gauss' triangular theorem. Since
$\{T_x : x\in\mathbb{N} \cup \{0\} \}=\{T_x: x\in\mathbb{Z}\}=\{x(2x+1) : x\in\mathbb{Z}\}$, 
 Gauss'  triangular theorem implies that the equation
$x(2x+1)+y(2y+1)+z(2z+1)=n$ has an integer solution $x,y,z$ for any non negative integer $n$. Inspired by this motivation, he defined in \cite{Sun}  that  for positive integers $a\leq b\leq c$,  
a ternary sum $\Phi_{a,b,c}(x,y,z)=x(ax+1)+y(by+1)+z(cz+1)$ is  {\it universal}  if the equation  
$$
x(ax+1)+y(by+1)+z(cz+1)=n
$$
has an integer solution for any non negative integer $n$.  He showed that  if $\Phi_{a,b,c}(x,y,z)$ is universal, then $(a,b,c)$ is one of the following $17$ triples:
$$ 
\begin{array}{llllll}
(1,1,2),&(1,2,2),&(1,2,3),&(1,2,4),&(1,2,5),&(2,2,2),\\
(2,2,3),&(2,2,4),&(2,2,5),&(2,2,6),&(2,3,3),&(2,3,4),\\
(2,3,5),&(2,3,7),&(2,3,8),&(2,3,9),&(2,3,10)
\end{array}
$$
and proved the universalities of some candidates.  In fact,  the universality of each ternary sum of triangular numbers was proved except the following $6$ triples:
\begin{equation}\label{candi}
\begin{array}{llllll}
\textbf{(2,2,6)},&\textbf{(2,3,5)},&\textbf{(2,3,7)},&(2,3,8),&(2,3,9),&(2,3,10).
\end{array}
\end{equation}
He also conjectured that for each of  these six triples, $\Phi_{a,b,c}(x,y,z)$ is universal. 
 
In this article, we prove that if $(a,b,c)$ is one of the triples written in boldface in the above candidates, then  $\Phi_{a,b,c}(x,y,z)$  is, in fact, universal.
To prove the universality, we use the method developed in  \cite{pentagonal}.  We briefly review this method in Section 2 for those who are  unfamiliar with it.  


Let $f$ be a positive definite integral ternary quadratic form. A symmetric matrix corresponding to $f$ is denoted by $M_f$. For a diagonal form $f(x,y,z)=a_1x^2+a_2y^2+a_3z^2$, we simply write $M_f=\langle a_1,a_2,a_3 \rangle$.
The genus of $f$ is the set of all integral forms that are locally isometric to $f$, which is denoted by $\text{gen}(f)$.  The set of all integers that are represented by the genus of $f$ ($f$ itself) is denoted by $Q(\text{gen}(f))$ ($Q(f)$, respectively).  For an integer $a$, we define 
 $$
 R(a,f)=\{(x,y,z)\in \mathbb{Z}^3 : f(x,y,z)=a\} \quad \text{and} \quad r(a,f)=|R(a,f)|.
 $$
  For a set $S$, we define $\pm S=\{ s : s \in S \ \text{or}\ -s \in S\}$. For two $v=(v_1,v_2,v_3), v'=(v'_1,v'_2,v'_3) \in  \z^3$ and a positive integer $s$,  we write $v\equiv v' \pmod s$ if $v_i\equiv v'_i \pmod s$ for any $i\in\{1,2,3\}$.

Any  unexplained notations and terminologies can be found in \cite{ki} or \cite{om}.

\section{General tools}
Let $a\leq b\leq c$ be positive integers. Recall that a ternary sum 
$$
\Phi_{a,b,c}(x,y,z)=x(ax+1)+y(by+1)+z(cz+1)
$$
is said to be universal if the equation $\Phi_{a,b,c}(x,y,z)=n$ has an integer solution for any non negative integer $n$. Note that  $\Phi_{a,b,c}(x,y,z)$ is universal if and only if the equation
$$
bc(2ax+1)^2+ac(2by+1)^2+ab(2cz+1)^2=4abcn+bc+ac+ab =: \Psi_{a,b,c}(n)
$$
has an integer solution for any non negative integer $n$. This is equivalent to the existence of an integer solution 
$X,Y$ and $Z$ of the  diophantine equation
$$
bcX^2+acY^2+abZ^2=\Psi_{a,b,c}(n)
$$
satisfying the following congruence condition
$$
X\equiv1 \ (\text{mod } 2a),~ Y\equiv1 \ (\text{mod } 2b) \text{ and } Z\equiv1 \ (\text{mod }2c).
$$

In some particular cases, representations of quadratic forms with some congruence condition corresponds to representations of a subform which is  suitably taken (for details, see \cite{poly}).

For the representation of an arithmetic progression by a ternary quadratic form,  
 we use the method developed in \cite{poly}, \cite{regular} and \cite{pentagonal}.   We briefly introduce this method for those who are unfamiliar with it. 
 
 Let $d$ be a positive integer and let $a$ be a non negative integer $(a\leq d)$. We define 
$$
S_{d,a}=\{dn+a \mid n \in \mathbb N \cup \{0\}\}.
$$
For integral ternary quadratic forms $f,g$, we define
$$
R(g,d,a)=\{v \in (\z/d\z)^3 \mid vM_gv^t\equiv a \ (\text{mod }d) \}
$$
and
$$
R(f,g,d)=\{T\in M_3(\mathbb{Z}) \mid  T^tM_fT=d^2M_g \}.
$$
A coset (or, a vector in the coset) $v \in R(g,d,a)$ is said to be {\it good} with respect to $f,g,d \text{ and }a$ if there is a $T\in R(f,g,d)$ such that $\frac1d \cdot vT^t \in \z^3$.  The set of all good vectors in $R(g,d,a)$ is denoted by $R_f(g,d,a)$.
If every vector contained in $R(g,d,a)$ is good, we write  $g\prec_{d,a} f$. If $g\prec_{d,a} f$, then by Lemma 2.2 of \cite{regular},  we have 
$$
S_{d,a}\cap Q(g) \subset Q(f).
$$
Note that the converse is not true in general.

\begin{thm}\label{maintool}
Under the same notations given above, assume that there is a partition $R(g,d,a)-R_f(g,d,a)=(P_1\cup \cdots \cup P_k) \cup (\widetilde{P}_1\cup \cdots \cup \widetilde{P}_{k'})$ satisfying the following properties: for each $P_i$, there is a  $T_i\in M_3(\mathbb Z)$ such that
\begin{enumerate}
\item[(i)] $T_i$ has an infinite order;
\item[(ii)] $T_i^tM_gT_i=d^2M_g$;
\item[(iii)]  for any vector $v \in \mathbb Z^3$ such that  $v\, (\text{mod }d)\in P_i$, $\frac1d\cdot vT_i^t\in \mathbb Z^3$ and $\frac1d\cdot vT_i^t\pmod d \in P_i \cup R_f(g,d,a)$, 
\end{enumerate}
and for each $\widetilde{P}_j$, there is a $\widetilde{T}_{j}\in M_3(\mathbb Z)$ such that
\begin{enumerate}
\item[(iv)] $\widetilde{T}_{j}^t M_g \widetilde{T}_{j}=d^2 M_g$;
\item[(v)]  for any vector $v \in \mathbb Z^3$  such that $v\, (\text{mod }d) \in \widetilde{P}_j$,  $\frac1d\cdot v\widetilde{T}_{j}^t \in \z^3$ and  $\frac1d\cdot v\widetilde{T}_{j}^t \pmod d \in P_1\cup\cdots\cup P_k \cup R_f(g,d,a)$. 
\end{enumerate}
Then we have
$$
(S_{d,a} \cap Q(g)) - \bigcup_{i=1}^{k} g(z_i)\mathcal S\subset Q(f),
$$
where the vector $z_i$ is a primitive eigenvector of $T_i$ and $\mathcal S$ is the set of squares of integers. 
\end{thm}

\begin{proof}
See Theorem 2.3 in \cite{poly}.
\end{proof}

\section{Universality of $x(ax+1)+y(by+1)+z(cz+1)$}
In this section, we prove that if $(a,b,c)=(2,2,6), (2,3,5)$ or $(2,3,7)$, then $\Phi_{a,b,c}(x,y,z)$ is universal. For the list of all candidates, see \eqref{candi} in the introduction.

\begin{lem} \label{case1-lem}  For any positive integer $n$, there are integers $a,b$ and $c$ such that $a^2+3b^2+3c^2=24n+7$ and $a\equiv b \equiv c \pmod 4$, whereas 
$$
a \not \equiv b \pmod 8 \quad \text{or} \quad a \not \equiv c \pmod8.
$$ 
\end{lem}

\begin{proof}  Since the class number of $f(x,y,z)=x^2+3(x-4y)^2+3(x-4z)^2$ is one and $S_{24,7}\subset Q(\text{gen}(f))$,  there is an integer solution $x^2+3y^2+3z^2=24n+7$ with $x\equiv y\equiv z\pmod 4$.  
If we define $g(x,y,z)=x^2+(x-8y)^2+3(x-8z)^2$, then it suffices to show that for any positive integer $n$,
$$
r(24n+7,f)-r(24n+7,g)>0.
$$ 
Note that  the genus of $g$ consists of
$$
M_{g}=\begin{pmatrix} 7&-1&3\\-1&55&27\\3&27&111\end{pmatrix}, \  M_2=\begin{pmatrix} 15&6&3\\6&28&14\\3&14&103\end{pmatrix} \ \ \text{and}  \ 
 \  M_3=\begin{pmatrix} 28&10&2\\10&31&-13\\2&-13&55\end{pmatrix}. 
$$
Furthermore, one may easily show by using the Minkowski-Siegel formula that  
$$
r(24n+7,f)=r(24n+7,\text{gen}(f))=4r(24n+7,\text{gen}(g)),
$$
where
$$
r(24n+7,\text{gen}(g))=\frac14r(24n+7,g)+\frac14r(24n+7,M_2)+\frac12r(24n+7,M_3).
$$
Therefore we have 
$$
r(24n+7,f)-r(24n+7,g)=r(24n+7,M_2)+2r(24n+7,M_3).
$$
Hence it suffices to show that every integer of the form $24n+7$ that is represented by $g$ is also represented by $M_2$ or $M_3$. 

One may show by a direct computation that 
$$
R(g,8,7)-R_{M_3} (g,8,7)=\pm\{ (1,0,0), (3,0,0)\}.
$$
Let $P_1=R(g,8,7)-R_{M_3}(g,8,7)$ and $T_1=\begin{pmatrix}8&-4&0\\0&-2&-12\\0&6&4\end{pmatrix}$.
Since $P_1$ and $T_1$ satisfy all conditions in Theorem \ref{maintool}, 
$$
(S_{24,7}\cap Q(g))-7\mathcal S \subset Q(M_3).
$$
Assume that $24n+7=7t^2$ for some positive integer $t$. Then $t$ has a prime divisor $p$ greater than $3$.  Since $g$ represents $7$,
$24n+7=7t^2$ is represented by $M_2$ or $M_3$  by Lemma 2.4 of \cite{poly}.  Therefore we have $r(24n+7,M_2)+2r(24n+7,M_3)>0$ for any positive integer $n$, 
which completes the proof. \end{proof}

\begin{thm} 
The ternary sum $x(2x+1)+y(2y+1)+z(6z+1)$ is universal. 
\end{thm}

\begin{proof}  It suffices to show that for any positive integer $n$,
$$
3(4x+1)^2+3(4y+1)^2+(12z+1)^2=24n+7
$$
has an integer solution $x,y,z$. 
By Lemma \ref{case1-lem}, there is an integer $a,b,c$ such that 
$$
a^2+3b^2+3c^2=24n+7, \quad  a\equiv b\equiv c \pmod 4  \quad   \text{and}    \quad  a \not \equiv b \pmod 8. 
$$
By changing signs suitably, we may assume that $a\equiv b \equiv c \equiv 1 \pmod 4$.  If $a \equiv 1 \pmod 3$, then everything is trivial. Assume that $a \equiv 2 \pmod 3$.  
In this case, note that 
$$
\left( \frac {a-3b}2\right)^2+3\left(\frac{a+b}2\right)^2+3c^2=a^2+3b^2+3c^2=24n+7
$$
and 
$$
 \frac {a-3b}2 \equiv 1 \pmod 3 \quad \text{and}\quad  \frac {a-3b}2 \equiv 1 \pmod 4.
 $$
 The theorem follows directly from this. 
\end{proof}

\begin{thm}
The ternary sum $x(2x+1)+y(3y+1)+z(5z+1)$ is universal.
\end{thm}

\begin{proof}
Let $f(x,y,z)=6x^2+10(2y+z)^2+15z^2$. We show that every integer of the form $40n+31$ for some positive integer $n$ is represented by $f$.  Note that the genus of $f$ consists of 
$$
M_f=\langle6\rangle\perp\begin{pmatrix}25&5\\5&25\end{pmatrix}, \ M_2=\langle 1,30,120 \rangle \  \  \text{and} \ \  M_3=\langle 10,15,24 \rangle.
$$
Note also that $S_{40,31}\subset Q(\text{gen}(f))$. By a direct computation, one may easily show  that $M_3\prec_{40,31}M_2$. Hence it suffices to show that every integer  of the form $40n+31$  that is represented by $M_2$ is also represented by $f$. 

Note that $R(M_2,40,31)-R_f(M_2,40,31)$ consists of the followings:
$$
\hspace{-3mm}
\begin{array}{lllllll}
\pm(1,\pm7,10), &\pm(1,\pm7,30), &\pm(1,\pm13,0), &\pm(1,\pm13,20), \\ 
\pm(9,\pm3,0), &\pm(9,\pm3,20), &\pm(9,\pm17,10), &\pm(9,\pm17,30), \\ 
\pm(11,\pm3,10), &\pm(11,\pm3,30), &\pm(11,\pm17,0), &\pm(11,\pm17,20), \\ 
\pm(19,\pm7,0), &\pm(19,\pm7,20), &\pm(19,\pm13,10), &\pm(19,\pm13,30). \\ 
\end{array}
$$
 Now, we define a partition of  $R(M_2,40,31)-R_f(M_2,40,31)$ as follows:
$$
\begin{array}{ll}
&P_1=\{(m_1,m_2,m_3) \mid 2m_1+m_2-m_3\equiv0\pmod5\},\\
&\widetilde{P}_1=\{(m_1,m_2,m_3) \mid 3m_1+m_2-m_3\equiv0\pmod5\},
\end{array}
$$
and we also define
$$
T_1=\begin{pmatrix}-20&-180&120\\-6&22&12\\-1&-3&-38\end{pmatrix}, \ \widetilde{T}_1=\begin{pmatrix}-20&180&-120\\-6&-22&-12\\-1&3&38\end{pmatrix}.
$$
Then one may easily show that it satisfies all conditions in Theorem \ref{maintool}.  Note that the eigenvector of $T_1$ is $(-3,1,0)$ and $120n+31$ is not of the form $39t^2$  for any positive integer $t$. Therefore for any non negative integer $n$,  the equation $f(x,y,z)=120n+31$ has an integer solution $(x,y,z)=(a,b,c)$.  Since $a\equiv 2b+c\equiv c\equiv 1 \pmod2, 2b+c\not\equiv0 \pmod3$ and $a\equiv\pm1\pmod5$, the equation
$$
15(4x+1)^2+10(6y+1)^2+6(10z+1)^2=120n+31
$$
has an integer solution for any non negative integer $n$. This completes the proof.
\end{proof}

\begin{thm}
The ternary sum $x(2x+1)+y(3y+1)+z(7z+1)$ is universal.
\end{thm}

\begin{proof}
Let $f(x,y,z)=6x^2+14y^2+21z^2$. We show that every integer of the form $21n+20$ for some positive integer $n$ is represented by $f$.
Note that  the genus of $f$ consists of  
$$
M_f=\langle 6,14,21 \rangle, \ M_2=\begin{pmatrix}5&1&2\\1&17&-8\\2&-8&26\end{pmatrix}, \ \ \text{and} \ \ M_3=\langle 3,14,42 \rangle.
$$
Note also that $S_{21,20}\subset Q(\text{gen}(f))$. 
One may easily show by a direct computation that $M_3\prec_{21,20}M_2$.  Hence it suffices to show that every integer of the form $21n+20$ that is represented by $M_2$ is also represented by $f$. 

One may show by a direct computation that the set $R(M_2,21,20)-R_f(M_2,21,20)$ consists of $52$ vectors, which is the union of the following sets:
$$
\begin{array}{lll}
&P_1=\pm\{(0,0,2), (0,0,5)\}, &P_2=\pm\{(0,4,0), (0,10,0)\},\\
&P_3=\pm\{ (3,2,0), (3,16,0)\}, &P_4=\pm\{ (2,0,0), (5,0,0)\},\\
& \\
&\widetilde{P}_1=\pm\{(6,9,10), (6,9,17)\}, &\widetilde{P}_2=\pm\{ (2,19,19),  (5,16,16)\},\\
&\widetilde{P}_3=\pm\{ (2,13,4),  (5,1,10)\}, &\widetilde{P}_4=\pm\{ (2,1,10),  (5,13,4)\},\\
&\widetilde{P}_5=\pm\{(3,7,12), (3,14,12)\}, &\widetilde{P}_6=\pm\{(6,5,12),(6,19,12)\},\\
&\widetilde{P}_7=\pm\{ (3,0,5), (3,0,19)\}, &\widetilde{P}_8=\pm\{(1,2,2), (8,16,16)\},\\
&\widetilde{P}_9=\pm\{(9,9,10), (9,9,17)\}.
\end{array}
$$
We also define
$$
{\setlength\arraycolsep{2pt}
\begin{array}{llll}
&T_1=\begin{pmatrix}-13&-24&-21\\-6&-3&21\\4&-12&21\end{pmatrix},
&T_2=\begin{pmatrix}-3&-40&16\\6&3&24\\-6&4&11\end{pmatrix},
&T_3=\begin{pmatrix}-9&30&-39\\12&9&3\\6&-6&-9\end{pmatrix},\\
&T_4=\begin{pmatrix}-1&36&0\\-12&-9&0\\-2&-12&21\end{pmatrix},
&T_5=\begin{pmatrix}-9&0&-45\\4&-21&13\\10&0&1\end{pmatrix},
&T_6=\begin{pmatrix}-1&-30&-15\\6&-9&27\\10&6&3\end{pmatrix},\\
&T_7=\begin{pmatrix}-9&24&12\\4&15&-24\\10&6&3\end{pmatrix},
&T_8=\begin{pmatrix}-3&36&-3\\6&-9&27\\-6&-12&15\end{pmatrix},
&T_9=\begin{pmatrix}-21&0&-21\\0&7&21\\0&-14&21\end{pmatrix}.
\end{array}}
$$
Note that  for any $i=1,2,\dots,9$, 
$$
T_i^t M_2 T_i =21^2 M_2 \quad  \text{and}  \quad (21^2\cdot T_i^{-1})^t M_2(21^2\cdot T_i^{-1})=21^2\cdot M_2.
$$  
Let $G=R_f(M_2,21,20)$, which is the set of good vectors. For any $i$, let  $Q_i$ be  the set $P_k$ or $\widetilde{P}_s$ for some $k$ or $s$. Let $U$ (and $V$) be an integral matrix $T_k$  or  $21^2T_s^{-1}$ for some $k$ or $s$. 
For any vector $v \in \z^3$ such that $v \, (\text{mod } 21) \in Q_1$, if 
$$
\frac1{21}\cdot vU^t \in \z^3\quad \text{and} \quad \frac1{21} \cdot vU^t  \, (\text{mod } 21)  \in G \cup  Q_{j_1} \cup \cdots \cup Q_{j_k},
$$
then we write 
$$
Q_1 \xrightarrow {U} G \cup  Q_{j_1} \cup \cdots \cup Q_{j_k}.
$$
We also use the notation
$$
Q_1 \xrightarrow {U} G\cup Q_2 \cup B  \xrightarrow {V}  G \cup  Q_{j_1} \cup \cdots \cup Q_{j_k},
$$
if  for any vector $v \in \z^3$ such that $v \, (\text{mod } 21) \in Q_1$, $\frac1{21}\cdot vU^t \in \z^3$ and either   $\frac1{21} \cdot vU^t  \, (\text{mod } 21)  \in G \cup Q_2$ or 
$$
\frac1 {21} \left(\frac1{21}\cdot vU^t\right)V^t \in \z^3 \quad \text{and} \quad  \frac1 {21} \left(\frac1{21}\cdot vU^t\right)V^t \in G \cup  Q_{j_1} \cup \cdots \cup Q_{j_k}.
$$
Here $B$ stands for the set of vectors of the form $\frac1{21} \cdot vU^t  \, (\text{mod } 21)$  that are not contained in  $G \cup Q_2$.

Now, one may show by direct computations that 
$$
\begin{array} {rll} 
& P_1 \xrightarrow{T_1} G \cup B \xrightarrow{T_2} G \cup \widetilde{P}_2,  &\widetilde{P}_2 \xrightarrow{T_3} G\cup B \xrightarrow{T_4} G \cup P_1, \\
& P_2 \xrightarrow{T_5} G \cup P_2 \cup B \xrightarrow{T_6} G,   & P_4 \xrightarrow{T_9} G \cup P_4, \\
 & \widetilde {P}_1 \xrightarrow{T_1} G \cup B \xrightarrow{T_2} G, & \widetilde{P}_3 \xrightarrow{21^2 T_1^{-1}} G \cup P_1 \cup \widetilde {P}_1,  \\
 & \widetilde{P}_4 \xrightarrow{21^2 T_2^{-1}} G \cup \widetilde{P}_2 \cup \widetilde {P}_3, &\widetilde{P}_5 \xrightarrow{T_5} G\cup P_2 \cup B \xrightarrow{T_6} G, \\
   & P_3 \xrightarrow{T_7} G \cup \widetilde{P}_3 \cup B \xrightarrow{T_8} G \cup P_3 \cup \widetilde{P}_5,  & \widetilde{P}_6 \xrightarrow{21^2 T_5^{-1}} G \cup P_2 \cup \widetilde {P}_5,  \\
   & \widetilde{P}_7 \xrightarrow{21^2 T_8^{-1}} G \cup P_3 \cup \widetilde {P}_6, & \widetilde{P}_8 \xrightarrow{21^2 T_9^{-1}} G  \cup \widetilde {P}_7, \\
    & \widetilde{P}_9 \xrightarrow{21^2 T_8^{-1}} G  \cup P_3 \cup \widetilde {P}_6. \\
 
  \end{array}
$$
As a sample of the above computations, assume that $v=(21x,21y,21z+2)$ for some integers $x,y,z$. Note that $v \, (\text{mod } 21) \in P_1$ and  
$$
w=\frac 1{21} \cdot vT_1^t=(-13 x-24 y-21 z-2,-6 x-3 y+21 z+2,4 x-12 y+21 z+2).
$$
 For any integers $x,y,z$, the vector $w \, (\text{mod } 21)$ is contained in $G$ or  
 $$
 w \, (\text{mod } 21)=(16,20,11) \quad  \text{or}  \quad  w \, (\text{mod } 21)=(19,2,2).
 $$
 If  $w \, (\text{mod } 21)=(16,20,11)$, then there are some integers $a,b$ such that $(x,y,z)=(18+3a,a,b)$. Hence $w=(-236-63a-21b, -106-21a+21b, 74+21b)$. Now, one may show by a direct
 computation that 
 $$
 \frac 1{21}\cdot wT_2^{t} \ (\text{mod } 21)  \in G \quad \text{or} \quad  \frac 1{21}\cdot wT_2^{t} \ (\text{mod } 21)=(19,2,2) \in \widetilde{P}_2.
 $$
If  $w \, (\text{mod } 21)=(19,2,2)$, then there are some integers $c,d$ such that $(x,y,z)=(3c,c,d)$. Hence $w=(-63c-21d-2, -21c+21d+2, 21d+2)$. Now, one may show by a direct
 computation that 
 $$
 \frac 1{21}\cdot wT_2^{t} \ (\text{mod } 21)  \in G \quad \text{or} \quad \frac 1{21}\cdot wT_2^{t} \ (\text{mod } 21)=(19,2,2)\in \widetilde{P}_2.
 $$

Note that the order of each $T_4T_3T_2T_1, T_5,T_8T_7$ and $T_9$ is infinite. Therefore, though the situation is slightly different from Theorem \ref{maintool}, we may still apply it to this case. The primitive integral eigenvectors of  these four matrices are $(0,0,1),(0,1,0),(3,2,0)$ and $(1,0,0)$, respectively. 
Therefore if $168n+41$ is not of the form $17t^2$ for some positive integer $t$, then $168n+41$ is represented by $f$. 
Assume that $168n+41=17t^2$ for some positive integer $t$. Then $t$ has a prime divisor relatively prime to $2\cdot3\cdot7$. 
Since $M_2,M_3\in\text{spn}(f)$, and both $M_2$ and $M_3$ represent $17$, $f$ represents $17t^2$ by Lemma 2.4 of \cite{poly}.  
Therefore the equation $f(x,y,z)=168n+41$ has an integer solution $a,b,c$ for any non negative integer $n$. 
Since $a\equiv b\equiv c \equiv 1\pmod 2$, $b\not\equiv0\pmod3$ and $a\equiv \pm1\pmod7$, the equation
$$
21(4x+1)^2+14(6y+1)^2+6(14z+1)^2=168n+41
$$ 
has an integer solution for any non negative integer $n$. This completes the proof.
\end{proof}

\end{document}